\documentclass[11pt]
{amsart}
\usepackage{amsthm,amsmath,amssymb,amsfonts}
\usepackage[colorlinks]{hyperref}
\usepackage[margin=1.4in]{geometry}
\usepackage{stmaryrd}
\usepackage{graphicx}
\newtheorem{theorem}{Theorem}
\newtheorem{thm}[theorem] {Theorem}

\newtheorem{lem}[theorem]{Lemma}

\newtheorem{cor}[theorem]{Corollary}

\newtheorem{obs}[theorem]{Observation}

\newtheorem{question}[theorem]{Question}
\newtheorem{example}[theorem]{Example}

\let\eps=\varepsilon

\let\polishlcross=\l
\def\l{\ifmmode\ell\else\polishlcross\fi}

\DeclareMathOperator{\Hom}{Hom}

\newcommand{\obsref}[1]{Observation~\ref{#1}}
\newcommand{\corref}[1]{Corollary~\ref{#1}}
\newcommand{\lemref}[1]{Lemma~\ref{#1}}
\newcommand{\exref}[1]{Example~\ref{#1}}
\newcommand{\thmref}[1]{Theorem~\ref{#1}}
\newcommand{\secref}[1]{Section~\ref{#1}}

\newcommand{\figref}[1]{Figure~\ref{#1}}

\newcommand{\ex}{\text{ex}}
\newcommand{\prob}{\ensuremath{\mathbb{P}}}

\newcommand{\rat}{\ensuremath{\mathbb{Q}}}
\newcommand{\real}{\ensuremath{\mathbb{R}}}
\newcommand{\nat}{\ensuremath{\mathbb{N}}}

\def\cA{{\mathcal A}}

\def\s{\sqrt{np}}

\def\cE{{\mathcal E}}

\def\cF{{\mathcal F}}

\def\A{{\mathcal A}}

\def\cK{{\mathcal K}}

\def\a{\alpha}

\def\d{\delta}

\def\f{\phi}

\def\g{\gamma}

\def\k{\kappa}

\def\l{\lambda}

\def\r{\rho}

\def\s{\sigma}

\newcommand{\proofstart}{{\bf Proof\hspace{2em}}}

\newcommand{\proofend}{\hspace*{\fill}\mbox{$\Box$}\vspace*{7mm}}

\newcommand{\rt}{\right}

\newcommand{\lt}{\left}

\newcommand{\hn}{\lt(H_n\rt)_{n\in \mathbb{N}}}

\usepackage{color}

\def\projectweb{\url{http://honza.ucw.cz/proj/erdos\_sos-k4e}}

\DeclareMathOperator{\im}{Im}

\def\Gerc{{
The work leading to this invention has received funding from the European Research Council under the European
Union's Seventh Framework Programme (FP7/2007-2013)/ERC grant agreement no.~259385.}}
\def\Ggauk{{GAUK 601812}}

\newdimen\margin
\def\textno#1&#2\par{
    \margin=\hsize
    \advance\margin by -4\parindent
           \setbox1=\hbox{\sl#1}
    \ifdim\wd1 < \margin
       $$\box1\eqno#2$$
    \else
       \bigbreak
       \hbox to \hsize{\indent$\vcenter{\advance\hsize by -3\parindent
       \sl\noindent#1}\hfil#2$}
       \bigbreak
    \fi}

\begin{document}

\title[A problem of Erd\H os and S\'os on $3$-graphs]{A problem of Erd\H os and S\'os on $3$-graphs}
\thanks{
{\Gerc}
{The first author was also supported by DFG within the research training group "Methods for Discrete Structures".}
{The third author was also supported by the student grant \Ggauk.}
}
\author[R.~Glebov]{Roman Glebov}
\address{
Department of Mathematics, ETH, 8092 Zurich, Switzerland.
Previous affiliations:
Mathematics Institute and DIMAP, University of Warwick, Coventry CV4 7AL, UK, and
Institut f\"{u}r Mathematik, Freie Universit\"at Berlin, Arnimallee 3-5, D-14195 Berlin, Germany.}
\email{roman.l.glebov@gmail.com}

\author[D.~Kr\'al']{Daniel Kr\'al'}
\address{
Mathematics Institute, DIMAP and Department of Computer Science, University of Warwick, Coventry CV4 7AL, UK. Previous affiliation:
Computer Science Institute, Faculty of Mathematics and Physics, Charles University, Prague, Czech Republic.}
\email{{\tt D.Kral@warwick.ac.uk}}

\author[J.~Volec]{Jan Volec}
\address{
Department of Mathematics, ETH, 8092 Zurich, Switzerland.
Previous affiliations:
Mathematics Institute and DIMAP, University of Warwick, Coventry CV4 7AL, and
Computer Science Institute, Faculty of Mathematics and Physics, Charles University, Prague, Czech Republic.}
\email{{\tt honza@ucw.cz}}

\begin{abstract}
We show that for every $\eps > 0$ there exist $\delta > 0$ and $n_0 \in \nat$ such that
every $3$-uniform hypergraph on $n\ge n_0$ vertices with the property that every $k$-vertex subset,
where $k \ge \delta n$, induces at least $\lt(\frac{1}{4} + \eps\rt)\binom{k}{3}$ edges,
contains $K_4^-$ as a subgraph, where $K_4^-$ is the $3$-uniform hypergraph on $4$ vertices
with $3$ edges. This question was originally raised by Erd\H os and S\'os.
The constant $1/4$ is the best possible.
\end{abstract}

\maketitle

\setcounter{footnote}{0}
\renewcommand{\thefootnote}{\fnsymbol{footnote}}

\section{Introduction}

One of the most influential results in the extremal graph theory is the celebrated theorem of Tur\'an~\cite{Tu41}
that determines the largest possible number of edges in an $n$-vertex graph without a complete subgraph of a~given size.
Erd\H os, Simonovits, and Stone~\cite{ErSi71, ErSt46} generalized Tur\'an's
theorem by showing that the {\em extremal number} of an arbitrary graph $F$, defined as
\[\ex(n,F)= \max\{|E(G)|:~G\mbox{ is a graph on } n \mbox{ vertices with no copy of } F\},\]
is asymptotically determined by its chromatic number. Specifically, for every graph $F$ with at least one edge,
\[\ex(n,F) = \lt(\frac{\chi(F)-2}{\chi(F)-1}+o(1)\rt)\binom{n}{2}.\]
Note that this problem is dual to determining the minimum number of edges $m(F)$
that guarantees that an $n$-vertex graph $G$ with at least $m(n,F)$ edges contains a copy of $F$,
since $m(n,F) = \ex(n,F) + 1$.
Simonovits~\cite{Simonovits68} showed the corresponding stability result, i.e.,
he proved that if a graph has almost the extremal number of edges and no copy of the forbidden graph $F$,
then it is close to the appropriate Tur\'an graph.
Specifically, it consists of $\chi(F)-1$ equal size nearly independent sets with almost all the edges between any pair of them.

In Ramsey-Tur\'an type problems, which were introduced by S\'os
in~\cite{Sos69}, we are again interested in the largest possible number of edges of
an $n$-vertex graph without a copy of $F$, but under the
additional restriction that the graph has to be somewhat ``far'' from the
Tur\'an graph. Typically, such a restriction is expressed by requiring that
the graph has only sublinear independence number. For more details on
Ramsey-Tur\'an theory, see the survey by Simonovits and
S\'os~\cite{SimonovitsSos01}.

Here, we are dealing with an even stronger restriction;
we require every linear-size subset of vertices not only to induce at least
one edge, but actually to induce a positive proportion of all possible edges.
We define the {\em $\delta$-linear density} of a graph $G$
to be the smallest density induced by a $\delta$-fraction of vertices, i.e.,
\[d(G,\delta)= \min\lt\{\frac{|E(G[A])|}{\binom{|A|}{2}}:~A\subseteq V(G), |A|\geq \delta |V(G)|\rt\}.
\]
Note that for any graph $G$ the function $d(G,\delta)$ is a non-decreasing function of $\delta$ taking values in $[0,1]$.

However, requiring a positive $\delta$-linear density immediately forces large graphs to behave similarly to random graphs
in the sense that they contain every given graph as a subgraph.
\begin{obs}
\label{obs:ldgraphs}
For every $\eps > 0$ and a fixed subgraph $F$, there exist $\delta > 0$ and $n_0 \in \nat$ such that
every graph $G$ on at least $n_0$ vertices with $d(G, \delta) \ge \eps$ contains $F$ as an subgraph.
\end{obs}
A proof of this observation follows, e.g., from~\cite[Theorem~1]{Rodl86}.

The notion of the $\delta$-linear density has a natural generalization to a $k$-uniform hypergraph $H$,
where we define
\[
d(H,\delta) = \min\lt\{\frac{|E(H_n[A])|}{\binom{|A|}{k}}:~A\subseteq V(H), |A|\geq \delta |V(H)|\rt\}.
\]
From now on, we mostly focus on $3$-uniform hypergraphs.
Let $K_4$ be the complete $3$-uniform hypergraph on $4$ vertices,
and $K_4^-$ the $3$-uniform hypergraph on $4$ vertices with $3$ edges.
Erd\H os and S\'os~\cite[Problem~5]{ErSo82} asked whether, analogously to~\obsref{obs:ldgraphs} for graphs,
a sufficiently large hypergraph with positive $\delta$-linear density contains a copy of $K_4$, or at least a copy of $K_4^-$.
However, F\"uredi observed that the following construction of Erd\H os and Hajnal~\cite{ErHa72} shows that the situation in $3$-uniform hypergraphs is completely
different (hence giving a negative answer to the above question).
\begin{example}
\label{ex:main}
Consider a random tournament $T_n$ on $n$ vertices. Let $H_n$ be the $3$-uniform hypergraph on the same vertex set
consisting of exactly those triples that span an oriented cycle in $T_n$.
\end{example}
One can check that in every hypergraph obtained in this way,
any four vertices span at most two edges, i.e., $H_n$ does not contain $K_4^-$ for every $n$.
On the other hand, for every $\delta>0$, the $\delta$-linear density of $H_n$
tends to $1/4$ when $n$ goes to infinity.
For the exact reference to the observation, we refer the reader to~\cite{FrRo88},
where Frankl and R\"odl cite a personal communication with F\"uredi in 1983.
Additional information about this problem and its history can be found in~\cite[Section~5]{SimonovitsSos01}.
It is also worth noting that in~\cite{Rodl86}, R\"odl presented a sequence of
$K_4$-free (but not $K_4^-$-free) $3$-uniform hypergraphs with $\delta$-linear density tending to $1/2$.

Our main result is the following theorem stating that the bound
$1/4$ on $\delta$-linear density for $K_4^-$-free sequences given
by~\exref{ex:main} is in fact the best possible.
This answers a question of Erd\H os from~\cite{Er90},
which is based on the original problem of Erd\H{o}s and S\'os, positively.
\begin{thm}
\label{thm:main}
For every $\eps > 0$ there exist $\delta > 0$ and $n_0 \in \nat$ such that
every $3$-uniform hypergraph $H$ on at least $n_0$ vertices with $d(H, \delta) \ge 1/4+\eps$
contains $K_4^-$ as a subgraph.
\end{thm}

Our proof is based on the flag algebra calculus. However, unlike in the other proofs in the area,
we apply the calculus to get one of the intermediate results only. We have to apply additional
considerations to derive the desired statement. Further details can be found in Section 3, where
we give the overview of the proof.

\section{Flag Algebras}
The main tool used for the proof of~\thmref{thm:main} is the framework of flag algebras.
They were introduced by Razborov~\cite{Ra07} as a general tool to approach problems from extremal combinatorics.
The work of Razborov was inspired by the theory of dense graph limits of Borgs \emph{et
al.}~\cite{lovasz2, lovasz3} and Lov\'asz and
Szegedy~\cite{LoSz06}. However, unlike for the more powerful Lov\'asz-Szegedy theory,
in order to create the theory of flag algebras for dense graphs, Razborov was able to avoid using any kind of regularity result.
In the case of $3$-uniform hypergraphs, flag algebras have already been successfully applied to
various problems. These include progress on the celebrated Tur\'an conjecture
about the minimum edge-density that forces containing $K_4$ as a~subhypergraph~\cite{Ra-tur10,Oleg-tur11,Ra-tur11,Ra-tur12},
several Tur\'an-type problems for $3$-uniform hypergraphs~\cite{BaTa12,flagmatic,flagmatic2},
and the so-called jumps in hypergraphs~\cite{BaTa11}.

Let us now introduce the terminology related to flag algebras
needed in this paper. Since we deal only with $3$-uniform hypergraphs
without a subhypergraph isomorphic to $K_4^-$, we restrict our attention only
to this particular case.
The central notions we are going to introduce are an algebra $\A$ and algebras
$\A^{\sigma}$, where $\sigma$ is a fixed $3$-uniform hypergraph not containing $K_4^-$.
In order to precisely describe these algebras, we first need to introduce
some additional notation.

Let $\cF$ be the set of all non-isomorphic $3$-uniform hypergraphs
that does not contain $K_4^-$. Next, for every $\ell\in\nat$, let $\cF_\ell\subset \cF$
be the set of $\ell$-vertex hypergraphs from $\cF$.
We say that $\ell$ is the {\em size} of a hypergraph $F\in\cF_\ell$ and denote it by $|F|$.
For $H\in\cF_\ell$ and $H'\in\cF_{\ell'}$, we set
$p(H,H')$ to be the probability that a randomly chosen subset of
$\ell$ vertices in $H'$ induces a subhypergraph isomorphic to $H$.
Note that $p(H,H')=0$ for $\ell' < \ell$.
Let $\real\cF$ be the set of all formal linear combinations
of elements of $\cF$ with real coefficients. Furthermore, let $\cK$ be the linear
subspace of $\real\cF$ generated by all linear combinations of the form
\[H-\sum_{H'\in\cF_{|H|+1}}p(H,H')\cdot H'.\]
Finally, we define $\cA$ to be the space $\real\cF$ factorized by $\cK$.

The space $\cA$ has naturally defined operations of addition and multiplication by a scalar (i.e., by a real number).
We now introduce multiplication inside $\cA$.
We start with defining it
on the elements of $\cF$ in the following way. For $H_1, H_2 \in \cF$, and $H\in\cF_{|H_1|+|H_2|}$,
we define $p(H_1, H_2; H)$ to be the probability that a randomly chosen subset of $V(H)$
of size $|H_1|$ and its complement induce in $H$ subhypergraphs isomorphic
to $H_1$ and $H_2$, respectively.
We set
\[H_1 \times H_2 = \sum_{H\in\cF_{|H_1|+|H_2|}}p(H_1,H_2;H) \cdot H.\]
The multiplication on $\cF$ has clearly the unique linear extension to $\real\cF$, which gives rise
to a well-defined multiplication also in the factor algebra $\cA$. A formal proof
of this can be found in~\cite[Lemma 2.4]{Ra07}.

Let us now move to the definition of an algebra $\cA^\sigma$, where $\sigma \in \cF$ is a fixed labelled hypergraph.
It follows basically the same lines as the definition of $\cA$.
Let $\cF^{\sigma}$ be the
set of $3$-uniform $K_4^-$-free hypergraphs $H$ with a fixed {\em embedding} of
$\sigma$, i.e., an injective mapping $\theta$ from $V(\sigma)$ to $V(H)$ such that
$\im(\theta)$ induces in $H$ a subhypergraph isomorphic to $\sigma$.  The elements of
$\cF^{\sigma}$ are usually called \emph{$\sigma$-flags} within the flag
algebra framework, $\im(\theta)$ are referred to as
the {\em labelled vertices} of the $\sigma$-flag, and the induced subhypergraph by the labelled vertices is called the \emph{type}.
Again, for every $\ell\in\nat$, we define
$\cF^{\sigma}_\ell\subset \cF^{\sigma}$ to be the set of the $\s$-flags from
$\cF^{\sigma}$ that have size $\ell$ (i.e., the $\s$-flags with the underlying hypergraphs having $\ell$ vertices).
Analogously to the
case for $\cA$, for $H, H' \in\cF^{\sigma}$,
we set $p(H,H')$ to be the
following probability:
a randomly chosen subset of $|H|-|\sigma|$ vertices in
$V(H')\setminus\theta'(V(\sigma))$ together with $\theta'(V(\sigma))$ induces a
subhypergraph that is isomorphic to $H$ through an isomorphism $f$ that preserves the embedding of $\sigma$.
In other words, $f(\theta') = \theta$.
Let $\real\cF^{\sigma}$ be the set of all formal
linear combinations of elements of $\cF^\sigma$ with real coefficients, and
$\cK^\sigma$ be the linear subspace of $\real\cF^\sigma$ generated by all
the combinations of the form \[H-\sum_{H'\in\cF^\sigma_{|H|+1}}p(H,H')\cdot H'.\]
We define $\cA^\sigma$ to be $\real\cF^\sigma$ factorised by $\cK^\sigma$.
We now describe the multiplication on $\cA^\sigma$.
For $H_1, H_2\in \cF^\s$ and $H\in \cF^\s_{|H_1+|H_2|-|\s|}$, similarly to the multiplication for $\cA$, we define
$p(H_1, H_2; H)$ to be the probability that a randomly chosen subset of $V(H)\setminus \theta(V(\sigma))$
of size $|H_1|-|\sigma|$ and a disjoint set from $V(H)\setminus \theta(V(\sigma))$ of size $|H_2|-|\sigma|$
extend $\theta(V(\sigma))$ in $H$ to subhypergraphs isomorphic to $H_1$ and $H_2$, respectively.
Again, by isomorphic here we mean that there is an isomorphism that preserves the embedding of $\sigma$.
This definition naturally extends to $\cA^\sigma$.

Consider an infinite increasing sequence $(H_n)_{n\in\nat}$ of hypergraphs from $\cF$.
The sequence $(H_n)_{n\in\nat}$ is \emph{convergent} if the probability $p(G,H_n)$
has a limit for every $G\in\cF$.
A standard compactness argument (e.g., using Tychonoff's theorem~\cite{Tyc30})
yields that every such infinite sequence has a convergent subsequence.
All the following results presented in the remainder of this section
can be found in~\cite{Ra07}.
Fix a convergent increasing sequence $(H_n)_{n\in\nat}$ of hypergraphs.
For every $G\in\cF$, we set $\phi(G) = \lim_{n\to\infty} p(G,H_n)$
and linearly extend $\phi$ to $\cA$.
We usually refer to the mapping as to the limit of the sequence.
The obtained mapping $\phi$ is a homomorphism from $\cA$ to $\real$.
Moreover, for every $F\in \cF$, we obtain $\phi(F)\geq 0$.
Let $\Hom^+(\cA, \real)$ be the set of all homomorphisms $\phi$ between the algebra $\cA$ and $\mathbb{R}$
that satisfy $\phi(F)\ge0$ for every $F \in \cF$.
It is interesting to see that indeed this set is the set of all limits of convergent sequences from $\cF$.

For $\sigma\in\cF$ and an embedding $\theta$ of $\sigma$ in $H_n$,
define $p^\theta_n(G)=p(G,H_n)$ for every $G\in\cF^\sigma$.
Picking $\theta$ at random gives rise to a probability distribution on mappings
from $\A^{\sigma}$ to $\real$, for each $n\in\nat$.
Since $p(G,H_n)$ converges (as $n$ tends to infinity) for every $G\in\cF^\sigma$,
the sequence of these probability distributions on mappings from $\A^{\sigma}$ to $\real$ also converges.
In fact, $\phi$ itself fully determines
the random distributions of the mappings $\phi^\sigma$ (for $\sigma$ with $\phi(\sigma) > 0$).
In what follows,
$\phi^\sigma$ will be a randomly chosen mapping from $\cA^{\sigma}$ to $\real$
based on this limit distribution. Any mapping $\phi^\sigma$ from the support
of the limit distribution is a homomorphism from $\cA^{\sigma}$ to $\real$.

The last notion we introduce is the \emph{averaging} (or downward) operator
$\llbracket\cdot\rrbracket_{\sigma}: \cA^{\sigma} \to \cA$.
For $G \in \cF^\sigma$, let $G^\emptyset$ be the (unlabeled) hypergraph in $\cF$ corresponding to $G$.
The operator $\llbracket\cdot\rrbracket_{\sigma}$ is the linear operator
defined on the elements $G \in \cF^\sigma$ by $\llbracket{G}\rrbracket_{\sigma} = p_G \cdot G^\emptyset$, where
$p_G$ is the probability that a random injective mapping
from $V(\sigma)$ to $V(G^\emptyset)$ is an embedding of $\sigma$ in $G^\emptyset$ yielding $G$.
The key relation between $\phi$ and $\phi^\sigma$ is the following:
\[
\forall G\in\cA^\sigma,\quad \phi\lt(\llbracket{G}\rrbracket_{\sigma}\rt)= \phi\lt(\llbracket{\sigma}\rrbracket_{\sigma}\rt) \cdot \int \phi^\sigma(G),
\]
where the integration is over the probability space given
by the limit random distribution on $\phi^\sigma$.
Therefore, if $\phi^\sigma(G)\ge 0$ almost surely,
then $\phi\lt(\llbracket{G}\rrbracket_{\sigma}\rt)\ge 0$.
In particular,
\begin{equation}
\forall G\in\cA^\sigma,\quad \phi\lt(\llbracket{G^2}\rrbracket_{\sigma}\rt)\ge 0.
\label{eq:star}
\end{equation}

\section{Structure of the paper and sketch of the proof}
In this section, we give a quick overview over the rest of the paper and we also introduce the remaining notation.
Since we always deal with sequences of hypergraphs, it seems natural to define
the {\em lower density} of an increasing sequence of hypergraphs $(H_n)_{n\in
\mathbb{N}}$ to be the smallest induced density of a linear-size subset
and denote it by
\[\l\lt((H_n)_{n\in\nat}\rt):=\lim_{\eps\rightarrow 0}\lt\{\liminf_{n\rightarrow \infty}d(H_n,\eps)\rt\}.\]
A sequence $\hn$ of the random hypergraphs from
\exref{ex:main} has lower density $1/4$ with probability one.

\secref{sec:proof} is devoted to the proof of \thmref{thm:main},
i.e., we show that the constant $1/4$ is best possible.
Our proof is partially computer assisted. However, it is worth pointing out here that we do not apply the Semidefinite method straightforwardly,
since we essentially need to include the specific restricted structure of the considered graphs with positive lower density.
We do so by first obtaining inequalities in the flag algebra language that hold for every potential counterexample.
These inequalities are then used as part of the linear combination to prove that the density of some specific subhypergraph is zero.
And this allows us to conclude a contradiction.

The proof is split into the following three steps:
\begin{itemize}
\item
In Step 1, we deduce some structural properties of a sequence of $3$-uniform hypergraphs
with lower density at least $1/4$. Namely, we prove that a specific collection of inequalities
is satisfied by any such sequence. Two particular inequalities from this collection
will be used then in Step 2.

\item
In Step 2, we use the semidefinite method from~\cite{Ra07} together with the two inequalities derived in Step 1
to formulate a specific semidefinite program. We used a computer software, namely an SDP library
CSDP~\cite{csdp}, to find an approximate solution of this program. The approximate solution
given by CSDP was then turned into an exact one by a careful rounding.
This was also done with computer assistance, using a mathematical software package Sage~\cite{sage}.
Since the size of the solution is quite large, we cannot present it here.
However, we made it available online at~\projectweb.

The rounded solution provides a proof that the induced densities of some specific flags of size $7$,
which we fully describe later, tend to zero in every sequence
with lower density at least $1/4$.
In this step, the flag algebra approach is essential.

\item
In Step 3, we first observe that, in particular, the main claim of Step 2 implies that the induced density
of one specific flag of size $5$ depicted in~\figref{fig-butterfly} tends to zero.
\begin{figure}
\begin{center}
\includegraphics[width=60mm]{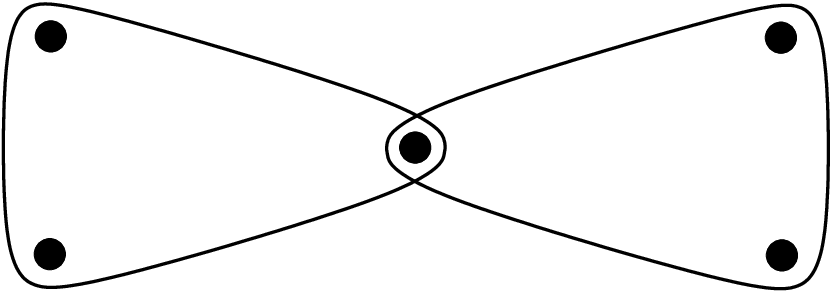}
\end{center}
\caption{The butterfly}
\label{fig-butterfly}
\end{figure}
Now suppose that~\thmref{thm:main} is false.
Hence, there exists a sequence of $3$-uniform hypergraphs without $K_4^-$ and with lower density larger than $1/4$.
We sparsify each element of the sequence uniformly to obtain a sequence of hypergraphs with lower density exactly $1/4$.
Therefore, the induced density of the butterfly in the sparsified sequence still tends to zero.
But this implies that the non-induced density of the butterfly in the original sequence tends to zero as well.
However, from that claim it will immediately follow that the edge density of the original sequence tends to zero, which is impossible.

\end{itemize}
Finally, in \secref{sec:concl} we give some concluding remarks, leaving some open questions.

\section{Proof of~\thmref{thm:main}}\label{sec:proof}

Suppose to the contrary that there exists $d>1/4$ and a sequence $\lt(H_n\rt)_{n\in \nat}$ of $3$-uniform $K_4^-$-free hypergraphs of increasing sizes
with lower density $d$.
Let us assume without loss of generality that $\lt(H_n\rt)_{n\in \nat}$ converges and $\f_C$ is the limit;
We aim to show that the edge density of $\f_C$ is zero,
contradicting the fact that it must be at least the lower density of the sequence.

As already described in the introduction, the proof is split into three steps.
\subsection{Deriving the inequalities}
In this subsection, we derive a collection of inequalities that are valid for any convergent sequence of 3-uniform hypergraphs $(H_n)_{n\in \mathbb{N}}$
with lower density at least $1/4$.

Let us first informally explain the main idea behind the inequalities.
Consider an arbitrary type $\sigma$ that has a positive density in the limit,
and let $F$ be a $\sigma$-flag with exactly one unlabelled vertex. Furthermore,
assume that $n$ is sufficiently large. Now fix a copy $S$ of $\sigma$ in $H_n$ and consider the set $U(S)$ of all vertices
$u \in V(H_n)$ that extends $S$ to a copy of $F$. The following two outcomes can happen:
\begin{itemize}
\item The size of $U(S)$ is small, i.e., $o(n)$. But then, for this choice of $S$,
both the probability that three random points belongs to $U(S)$,
and the probability that three random points belongs to $U(S)$ and span an edge, are $o(1)$.
\medskip
\item The size of $U(S)$ is large, i.e., $\Omega(n)$. But then, by the assumption on the lower density,
for this choice of $S$, the probability that three random points belong to $U(S)$ and span an edge
is asymptotically at least quarter of the probability that three random points belongs to $U(S)$.
\end{itemize}

Analysis of these two outcomes is given in the following lemma.
\begin{lem}
\label{lem:kappas}
Let $\lt(H_n\rt)_{n\in \nat}$ be a convergent sequence of 3-uniform hypergraphs with lower density at least $1/4$.
Let $\f$ be its limit, $\sigma$ a type, and $F \in \cF^\sigma_{v(\sigma)+1}$ a $\sigma$-flag.
In addition, let $\alpha_G$, for $G \in \cF^\sigma_{v(\sigma)+3}$, be the coefficients
such that
$F \times F \times F = \sum\limits_{G \in \cF^\sigma_{v(\sigma)+3}} \alpha_G \cdot G$.
Furthermore, let
\[
\kappa := \sum\limits_{G \in \cF^\sigma_{v(\sigma)+3}} \alpha_G \cdot G
\hfill\qquad\mbox{and}\qquad\hfill
\kappa^+ :=     \sum\limits_{G \in \cF^\sigma_{v(\sigma)+3}} \alpha^+_G \cdot G,
\]
where
\[
\alpha^+_G :=
\begin{cases}\alpha_G & \mbox{if the three non-labelled vertices of $G$ span an edge, and} \\0 & \mbox{otherwise.}\end{cases}
\]
It holds that
\begin{equation}\label{eq:kappas_unlab}
\phi\left(\llbracket 4\kappa^+-\kappa\rrbracket_\sigma\right)\geq 0.
\end{equation}
\end{lem}
\noindent Note that
$\alpha_G \in \{0,1\}$ for every $G\in\cF^\sigma_{v(\sigma)+3}$.

\proofstart
First observe that if $\phi\left(\sigma^\emptyset\right) = 0$, then all the
terms in the left-hand side of (\ref{eq:kappas_unlab}) are equal to zero.
For the rest of the proof, we assume $\phi\left(\sigma^\emptyset\right) > 0$.
Our aim is to show that following holds when $\phi^\sigma$ is drawn from the distribution given by $\phi$:
\begin{equation}\label{eq:kappas_lab}
\prob\left(\phi^\sigma(4\kappa^+-\kappa) \geq 0\right)=1.
\end{equation}
Since the inequality~\eqref{eq:kappas_unlab} is an immediate corollary
of~\eqref{eq:kappas_lab}, it is enough to establish~\eqref{eq:kappas_lab}.

Suppose for the contrary that~\eqref{eq:kappas_lab} is not true, i.e., $\prob\lt(\f^\sigma(4\k^+-\k) \leq -\eps\rt)\geq \eps$ for some $\eps>0$.
Then there exists $n_0 \in \nat$
such that
for every $n \ge n_0$ the graph $H_n$ contains a positive proportion of copies of $\sigma$
with the following property. For every such copy~$S_n$
of $\sigma$ in $H_n$, the set $U(S_n)=\{v\in V(H_n)\setminus V(S_n): H_n[V(S_n) \cup \{v\}] \cong F\}$
has size at least $\eps n$ and the number of edges of $H_n[U]$ is at most $\frac{1}{4}(1-\eps) \cdot \binom{|U|}{3}$.
But this contradicts the fact that $\l(H_n)\ge 1/4$.
\proofend

\begin{figure}
\begin{center}
\includegraphics[width=100mm]{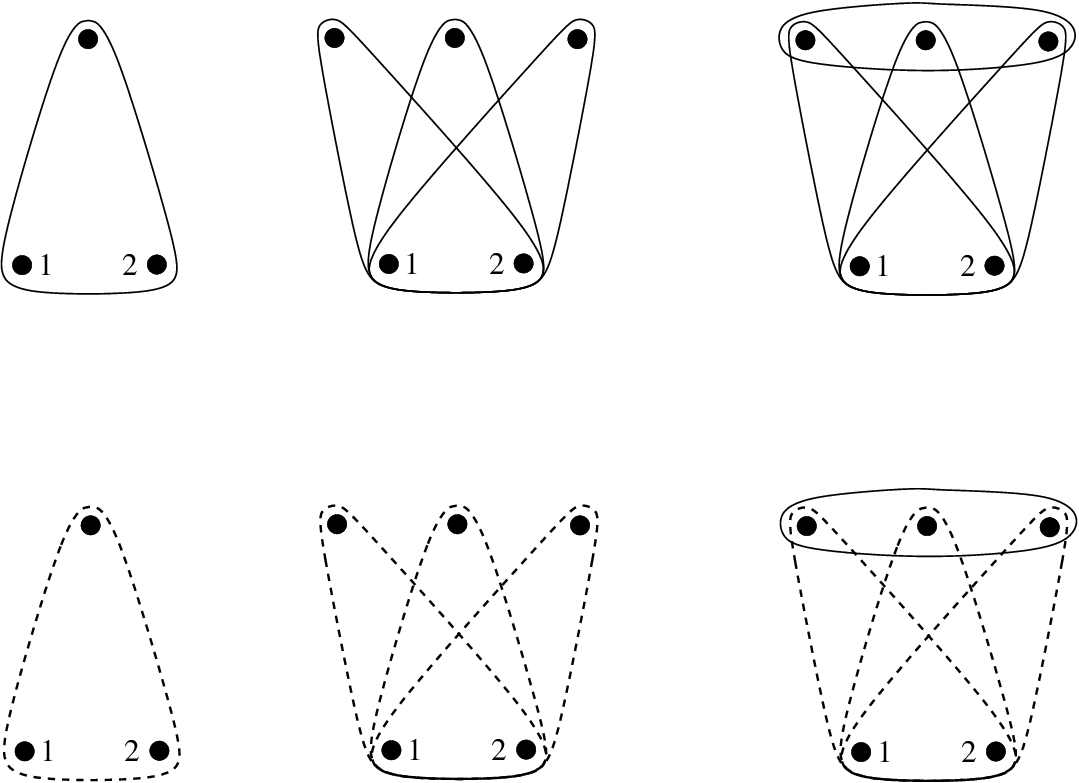}
\end{center}
\caption{The expressions used in inequalities (\ref{eq:k1}) and (\ref{eq:k2}): a solid line denotes a hyperedge, dotted line a non-hyperedge, and finally we sum over all the possible choices of triples
where there is neither a dotted nor a solid line.}
\label{fig-kappa}
\end{figure}

Now we are ready to present the two inequalities we are going to use in the next step.
Let~$\k_1$ be $\k$ and let $\k_1^+$ be $\k^+$
as in~\lemref{lem:kappas} for $F$ being the edge with two labelled vertices (hence $\sigma$ is the unique type on $2$ vertices;
from now on we write just $\bf 2$ when we use this type).
The first line in \figref{fig-kappa} shows this flag $F$ and the corresponding linear combinations $\k_1$ and $\k_1^+$.
Next, let~$\k_2$ be $\k$ and let $\k_2^+$ be $\k^+$
for $F$ being the non-edge with two labelled vertices (again, $\sigma$ is $\bf 2$).
The second line in \figref{fig-kappa} shows the corresponding flag $F$ and the combinations $\k_2$ and $\k_2^+$.

\lemref{lem:kappas} then implies that
\begin{equation}\label{eq:k1}
\f\lt(\llbracket 4\k_1^+-\k_1\rrbracket_{\bf 2}\rt)\geq 0,
\end{equation}
and
\begin{equation}\label{eq:k2}
\f\lt(\llbracket 4\k_2^+-\k_2\rrbracket_{\bf2}\rt)\geq 0.
\end{equation}

\subsection{Applying the Semidefinite Method}
Recall that $\cF_7$ is the set of all $3$-uniform $K_4^-$-free hypergraphs of size $7$.
It holds that $\lt|\cF_7\rt| = 8157$.
Let $\cE_7$ be the set of all hypergraphs from $\cF_7$
that can possibly appear as induced subgraphs of the hypergraphs from~\exref{ex:main}.
It holds that $\lt|\cE_7\rt| = 247$.

\begin{figure}
\begin{center}
\includegraphics[height=72mm]{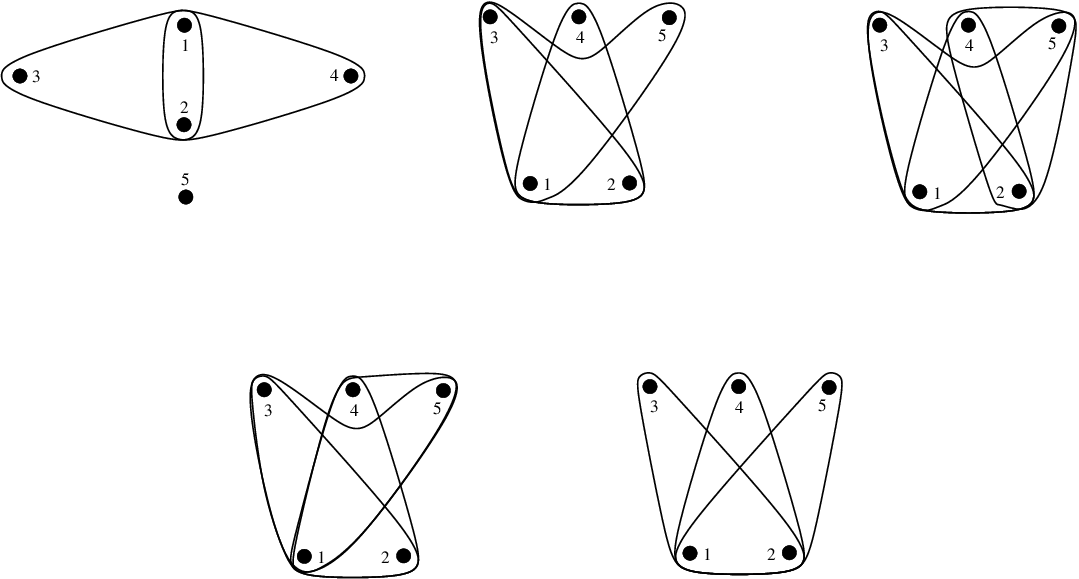}
\end{center}
\caption{The types $\s_4$ through $\s_8$: solid lines denote edges, and a triple without a solid line spans a non-edge.}
\label{fig-roots}
\end{figure}

The main lemma of this part of the proof is the following.
\begin{lem}
\label{lem:sdp}
There exist rational numbers $\lt(\a_G\rt)_{G\in \cF_7}$, where
$\a_G=0$ for $G \in \cE_7$, and $\a_G<0$ otherwise, such that the following holds.
If $\f$ is an element of $\Hom^+(\cA, \real)$
that satisfies~\eqref{eq:k1} and~\eqref{eq:k2},
then
\[\f\lt(\sum_{G \in \cF_7} \a_G \cdot G\rt) = 0\mbox{.}\]
\end{lem}
\proofstart
As we already mentioned, the proof was found with computer assistance.
An instance of an SDP program was used to find positive rationals $\g_1$ and $\g_2$ and $8$
symmetric positive semidefinite matrices $A_1,A_2,\dots A_8$ with rational entries
such that the following holds:
\[\f\lt(
\g_1 \cdot \llbracket 4\k_1^+-\k_1\rrbracket_{\bf2} +
\g_2 \cdot \llbracket 4\k_2^+-\k_2\rrbracket_{\bf2} +
\sum_{i=1,\dots,8} \llbracket x_i^T A_i x_i\rrbracket_{\sigma_i}\rt)
= \f\lt(\sum_{G \in \cF_7} \a_G \cdot G\rt), \]
where
\begin{itemize}
\item the type $\sigma_1$ is the type with one vertex,
\item the type $\sigma_2$ is the type with three vertices and no edge,
\item the type $\sigma_3$ is the type with three vertices and one edge,
\item the types $\sigma_4,\dots,\sigma_8$ are the five specific types on five vertices given in~\figref{fig-roots},
\item the vector $x_1 \in (\real\cF^{\sigma_1}_4)^{\lt|\cF^{\sigma_1}_4\rt|}$ is the vector whose
$j$-th coordinate is equal to the $j$-th element of the canonical base of $\real\cF^{\sigma_1}_4$,
\item for $i=2,3$, the vector $x_i \in (\real\cF^{\sigma_i}_5)^{\lt|\cF^{\sigma_i}_5\rt|}$
is the vector whose $j$-th coordinate is equal to the $j$-th element of the canonical base of $\real\cF^{\sigma_i}_5$,
\item for $i=4,5,\dots,8$, each vector $x_i \in (\real\cF^{\sigma_i}_6)^{\lt|\cF^{\sigma_i}_6\rt|}$
is the vector whose $j$-th coordinate is equal to the $j$-th element of the canonical base of $\real\cF^{\sigma_i}_6$, and
\item $\a_G=0$ if $G \in \cE_7$, and $\a_G<0$ otherwise.

\end{itemize}
Since the left-hand side of the above equation is non-negative by~\eqref{eq:star},~\eqref{eq:k1}, and~\eqref{eq:k2},
it follows that \[\f\lt(\sum_{G \in \cF_7} \a_G \cdot G\rt) \ge 0.\]
On the other hand, the values of $\alpha_G$ are non-positive for every $G \in \cF_7$, hence
the left-hand side of the last inequality is trivially also at most zero.
\proofend

To illustrate the complexity of the proof, we note that the sizes $a_i$ of the sets
$\cF^{\sigma_1}_4, \cF^{\sigma_2}_5, \cF^{\sigma_3}_5, \cF^{\sigma_4}_6, \cF^{\sigma_5}_6, \cF^{\sigma_6}_6, \cF^{\sigma_7}_6$, and $\cF^{\sigma_8}_6$,
which coincide with the order of the matrices $A_i$,
are $5, 95, 47, 191, 135, 95, 101$, and $148$, respectively.
The numerical values of $\g_1, \g_2$ and $A_i$'s can be downloaded from~\projectweb.

Each matrix $A_i$ is not stored directly but
as an appropriate number of vectors $v^j_i \in \rat^{a_i}$ and positive rationals
$w^j_i$ such that \[ A_i = \sum_{j=1}^{r_i} w^j_i \cdot v^j_i \times \lt(v^j_i\rt)^T,\]
where $r_1=3, r_2=45, r_3=21, r_4=63, r_5=60, r_6=43, r_7=31,$ and $r_8=25$.

In order to make an independent verification of our computations easier, we created a sage script called ``verify.sage'',
which is also available on the web page mentioned above. Instead of generating the lists of flags
and computing all the probabilities involved in our computations, the verification program
uses data from Flagmatic~1.5.1~\cite{flagmatic, flagmatic2} applied to
find an upper bound on the Tur\'an density of $K_4^-$.
However, we used these data only to simplify the verification part and they were not used
in our computations to find $\g_1, \g_2$ and $A_i$s (this was to decrease a chance of encountering a program-related error).

\lemref{lem:sdp} immediately yields the following.
\begin{cor}
\label{cor:sdp}
Let $F$ be a hypergraph of size at most $7$ such that its induced density in the hypergraphs from~\exref{ex:main} is zero,
and let $\f$ be as in~\lemref{lem:sdp}.
Then,
$\f(F)=0$.
\end{cor}

\subsection{The final step}
Let $B$ be the butterfly hypergraph depicted in~\figref{fig-butterfly}.
Since the hypergraphs from \exref{ex:main} are $B$-free, \corref{cor:sdp} implies that $\f(B)=0$
whenever a homomorphism $\f\in\Hom^+(\cA, \real)$ satisfies~\eqref{eq:k1} and~\eqref{eq:k2}. In particular, for the limit $\f$
of any convergent sequence of $K_4^-$-free hypergraphs with lower density at least $1/4$, we have $\f(B)=0$.

Recall that $\f_C$ is the limit of a hypothetical convergent sequence $\lt(H_n\rt)_{n\in \nat}$ that is a counterexample to~\thmref{thm:main}.
Therefore, $\l\lt(\lt(H_n\rt)_{n\in \nat}\rt) > 1/4$. Our aim is to show that the edge density of $\f_C$ is equal to zero,
which would contradict $\l\lt(\lt(H_n\rt)_{n\in \nat}\rt) > 1/4$.

Instead of applying the claim directly to $\f_C$, we
first construct from $\lt(H_n\rt)_{n\in \nat}$ a new sequence $\lt(H'_n\rt)_{n\in \nat}$
such that $\l\lt(\lt(H'_n\rt)_{n\in \nat}\rt) = 1/4$. This is done simply by
a \emph{random sparsification} of $\lt(H_n\rt)_{n\in \nat}$, i.e.,
removing each edge of each $H_n$ at random with the appropriately chosen probability.
\begin{obs}
\label{obs:ld(gp)}
Let $\lt(H_n\rt)_{n\in \nat}$ be a sequence of $3$-uniform $K_4^-$-free hypergraphs with $d := \l\lt(\lt(H_n\rt)_{n\in \nat}\rt) > 1/4$.
Furthermore, for every $n\in\nat$, let $H'_n$ be a random subhypergraph of $H_n$
obtained by removing every hyperedge of $H_n$ independently with probability $1-\frac{1}{4d}$. Then
\[\prob\lt(\l\lt(\lt(H'_n\rt)_{n\in \nat}\rt) = 1/4\rt)=1.\]
\end{obs}

Let $\lt(H'_n\rt)_{n\in \nat}$ be the sequence of hypergraphs
obtained from $\lt(H_n\rt)_{n\in \nat}$.
Let $\f'_C$ be its limit (the sequence $\lt(H'_n\rt)_{n\in \nat}$ must be convergent).
\lemref{lem:kappas} implies that $\f'_C$
satisfies both~\eqref{eq:k1} and~\eqref{eq:k2}. Hence, by \corref{cor:sdp},
the induced density of $B$ in $\phi'_C$ is equal to zero. 
But this implies that $\phi_C(F)=0$ for any $F \in \cF$ that contains
$B$ as a non-induced subgraph. This holds because otherwise there would be some $F\in \cF$ that contains
a non-induced copy of $B$ and $\phi_C(F) > 0$. However, a positive proportion of the induced copies of $F$
from $\phi_C$ will be then turned, by a random sparsification, to induced copies of $B$ in $\phi'_C$.
We conclude that the \emph{non-induced} density of the butterfly $3$-graph $B$
in $\phi_C$ is equal to zero. In particular, any non-negative
combination of the elements of $\cF_5$ that contain the butterfly as a
non-induced subgraph must be zero in $\phi_C$.

We are now ready to derive that the edge density of $\f_C$ is zero.
Let $\r \in \cF_3$ be the hypergraph on three (unlabelled) vertices that span an edge.
Let
$\r_1 \in \cF_3^1$ be the flag corresponding to $\rho$ with exactly one vertex labelled.

Observe that all elements of $\cF_5^1$ with positive coefficient in $\rho_1^2$ contain $B$ as a subgraph.
So, $\f_C\lt( \llbracket (\r_1)^2 \rrbracket_1 \rt) = 0$, which implies that
\[
\f_C\lt(\r\rt)^2 = \f_C\lt(\llbracket \r_1 \rrbracket_1 \rt)^2 \le
\f_C\lt( \llbracket (\r_1)^2 \rrbracket_1 \rt)=0.
\]

This completes the proof of \thmref{thm:main}.

\section{Concluding remarks and open questions}
\label{sec:concl}
In \thmref{thm:main}, we proved that for $3$-uniform hypergraphs, the threshold
lower density for every four vertices to induce at most two edges is $1/4$.
We have recently learned~\cite{Ro-priv} that one can prove the same result using
hypergraph regularity under the assumption that all parts in the partitions are
regular.

Recently, Falgas-Rarvy, Pikhurko and Vaughan~\cite{FalgPikVaug:K4-} used the
semidefinite method to prove
that the value $1/4$ is also the minimum
co-degree threshold for containing a copy of $K_4^-$. Furthermore, they showed
that the limit of the sequence of hypergraphs corresponding to cyclically
oriented triangles in random tournaments is in fact the unique $K_4^-$-free
homomorphism from $\Hom^+(\cA,\real)$, where the elements of the sequence have
minimum co-degree equal to $1/4+o(1)$. An adaptation of their approach to our
problem yields the same uniqueness result as they showed for the co-degree
condition.
\begin{theorem}
\label{thm:k4estab}
Let $(H_n)_{n\in\nat}$ be a convergent sequence of $K_4^-$-free $3$-graphs so that $\lambda((H_n)) = 1/4$,
and let $\phi \in \Hom^+(\cA,\real)$ be its limit.
It holds that $\phi$ is equal to the homomorphism $\psi$, which is the limit
of the sequence of $3$-graphs from Construction~\ref{ex:main} with increasing orders.
\end{theorem}
\noindent A proof of this theorem can be found in the Ph.D. thesis of the third author~\cite{volec:phd}.

Bhat and R\"odl~\cite{BhRo13+} applied our \thmref{thm:main} to show the following.
A number $\a>0$ is a {\em jump for $r$-uniform pseudorandom hypergraphs}, if there exists $\d>0$
such that every sequence of $r$-uniform pseudorandom hypergraphs of increasing sizes with limit density larger than $\a$
contains subgraphs of increasing sizes with density at least $\a+\d$.
Bhat and R\"odl showed that every number between $0$ and $0.3192$ is a jump for $3$-uniform pseudorandom hypergraphs.
Notice that, if Question~\ref{q:example} (given further) is answered positively for $r=4$,
then similarly to the proof of~\cite[Theorem~1.6]{BhRo13+} one can show that every number between $1/8=0.125$ and $0.1262869$ is a
jump for $4$-uniform pseudorandom hypergraphs.

A natural question is to seek an extension of our results for hypergraphs with larger uniformity.
Frankl and R\"odl~\cite{FrRo88} mention a construction of Frankl and F\"uredi ({\cite[Example~3]{FrFu84}})
as a possible source of examples of $r$-uniform hypergraphs with lower density $1/2^{r-1}$
with no $r+1$ vertices spanning three edges.
Here, we present a variant of this construction
that was also independently obtained by them~\cite{Ro-priv} which indeed gives such examples.

Every $(r-2)$-simplex in $\real^{r-1}$ has two possible {\em orientations} with respect to the hyperplane it spans.
We naturally call an $(r-1)$-simplex in $\real^{r-1}$ {\em compliant} if all its
$(r-2)$-dimensional faces are either all oriented towards the inside or all
towards the outside of the $(r-1)$-simplex.

\begin{example}
\label{ex:simplex}
Consider $n$ points in $\mathbb{R}^{r-1}$ in general convex position.  Orient
$(r-2)$-simplices spanned by these points uniformly at random (i.e., every of
the two possible orientations has probability $1/2$).  Let $r$ points form an
edge if the corresponding $(r-1)$-simplex is compliant.
\end{example}

The lower density of $r$-uniform hypergraphs from \exref{ex:simplex} tends to $1/2^{r-1}$.
The following argument shows that every $r+1$ vertices span at most two edges
in every such $r$-uniform hypergraph.
Let us fix $r+1$ vertices in $\real^{r-1}$ in general convex position and
consider an orientation of the $(r-2)$-simplices formed by them.
\begin{itemize}
\item
The intersection of two $(r-1)$-simplices is an $(r-2)$-simplex that is either part of the convex hull of the $r+1$ points
or not.
\item
Two compliant $(r-1)$-simplices whose intersection is part of the convex hull
have their $(r-2)$-dimensional faces
either all oriented towards their insides or all towards the outside.
\item
Similarly, out of two compliant $(r-1)$-simplices whose intersection is not part of the convex hull,
one has all its $(r-2)$-dimensional faces oriented towards its inside, and the other towards its outside.
\item
For every three $(r-1)$-simplices, there are two possibilities for the three
$(r-2)$-dimensional faces in their pairwise intersections:
either none of them or exactly two of them are in the convex hull.
So, not all three $(r-1)$-simplices can be compliant.
\end{itemize}
This leads us to the following question.
\begin{question}
\label{q:example}
Let $\hn$ be a convergent sequence of $K_{r+1}^3$-free $r$-uniform hypergraphs,
where $K_{r+1}^3$ is the $r$-uniform hypergraph on $r+1$ vertices with three edges.
Is it true that $\hn$ has lower density at most $1/2^{r-1}$?
\end{question}

\paragraph{\bf Acknowledgments}
The authors thank Oleg Pikhurko and Vojta R\"odl for fruitful discussions and
helpful comments on the results contained in this paper.
The authors also thank to the anonymous referee for carefully reading the
manuscript and for their valuable comments, which greatly improved the
presentation of the results.

\end{document}